\documentclass[reqno,12pt]{amsart}
\usepackage{amsfonts}
\usepackage{amssymb}
\usepackage{graphicx}
\usepackage{fullpage}
\usepackage{url}

\begin{document}

\newtheorem{theorem}{Theorem}
\newtheorem{lemma}{Lemma}
\newtheorem{proposition}{Proposition}
\newtheorem{remark}{Remark}


\def\eqref#1{(\ref{#1})}
\def\eqrefs#1#2{(\ref{#1}) and~(\ref{#2})}
\def\eqsref#1#2{(\ref{#1}) to~(\ref{#2})}
\def\sysref#1#2{(\ref{#1})--(\ref{#2})}

\def\Eqref#1{Eq.~(\ref{#1})}
\def\Eqrefs#1#2{Eqs.~(\ref{#1}) and~(\ref{#2})}
\def\Eqsref#1#2{Eqs.~(\ref{#1}) to~(\ref{#2})}
\def\Sysref#1#2{Eqs. (\ref{#1}),~(\ref{#2})}

\def\secref#1{Sec.~\ref{#1}}
\def\secrefs#1#2{Sec.~\ref{#1} and~\ref{#2}}

\def\appref#1{Appendix~\ref{#1}}

\def\Ref#1{Ref.~\cite{#1}}
\def\Refs#1{Refs.~\cite{#1}}

\def\Cite#1{${\mathstrut}^{\cite{#1}}$}

\def\EQ{\begin{equation}}
\def\endEQ{\end{equation}}


\def\fewquad{\qquad\qquad}
\def\severalquad{\qquad\fewquad}
\def\manyquad{\qquad\severalquad}
\def\manymanyquad{\manyquad\manyquad}

\def\downindex#1{{}_{#1}}
\def\upindex#1{{}^{#1}}

\def\mathtext#1{\hbox{\rm{#1}}}

\def\hp#1{\hphantom{#1}}

\def\parder#1#2{\partial{#1}/\partial{#2}}
\def\parderop#1{\partial/\partial{#1}}

\def\xder#1{{#1}\downindex{x}}
\def\vder#1{{#1}\downindex{v}}

\def\D#1{D\downindex{#1}}

\def\X#1{{\bf X}_{#1}}
\def\prX#1{{\bf X}_{#1}^{(1)}}
\def\Y#1{{\bf Y}_{#1}}
\def\sgn{{\rm sgn}}

\def\ie/{i.e.}
\def\eg/{e.g.}
\def\etc/{etc.}
\def\const{{\rm const.}}


\title{Linearizability Criteria for Systems of Two Second-Order
Differential Equations by Complex Methods}

\author{S. Ali}
\address{
School of Electrical Engineering and Computer Sciences,
National University of Sciences and Technology,
H-12 Campus, Islamabad 44000, Pakistan}
\email{sajid\_ali@mail.com}

\author{F. M. Mahomed}
\address{
School of Computational and Applied Mathematics, Centre for Differential Equations, Continuum Mechanics and Applications, University of the Witwatersrand, Wits 2050, South Africa }
\email{Fazal.Mahomed@wits.ac.za}

\author{Asghar Qadir}
\address{
Center For Advanced Mathematics and Physics, National
University of Sciences and Technology, Campus H-12, 44000,
Islamabad, Pakistan }
\email{aqadirmath@yahoo.com}

\begin{abstract}
Lie's linearizability criteria for scalar second-order ordinary
differential equations had been extended to systems of second-order
ordinary differential equations by using geometric methods. These
methods not only yield the linearizing transformations but also the
solutions of the nonlinear equations. Here complex methods for a
scalar ordinary differential equation are used for linearizing
systems of two second-order ordinary and partial differential
equations, which can use the power of the geometric method for
writing the solutions. Illustrative examples of mechanical systems
including the Lane-Emden type equations which have roots in the study of
stellar structures are presented and discussed.
\end{abstract}

\def\sep{;}
\keywords{
complex linearization, complex Lie symmetries, Lie linearizability
}

\maketitle

\section{ Introduction }
Lie developed a systematic procedure for solving nonlinear ordinary
differential equations (ODEs) with
some minimal symmetry under transformations of the dependent and
independent variables, called {\it point transformations}, by using
group theory \cite{lie1}. For second-order ODEs he provided criteria
for their being transformable to linear ODEs, provided they are
maximally symmetric \cite{lie2}. The requirement is that they be at
most cubically nonlinear in the first derivative and that their
coefficients satisfy a set of four constraints for the first
derivatives, involving two auxiliary functions. Tress\'{e}
\cite{tre1,tre2} reduced the number of constraining equations to two
for higher derivatives by eliminating the auxiliary functions.

Since Lie's time there have been various developments following up
his work. There was the development of contact transformations (see
e.g., \cite{ba,ibr}) and their use for linearizing two classes of
third-order ODEs \cite{gre1,gre2}; approach of Cartan type used for
the same purpose \cite{che1,che2}; a method for a special case of a
class of third-order scalar ODEs \cite{mel}; Lie's algebraic methods
used for classification of systems of ODEs \cite{ml1,ml2,wm};
extension to partial differential equations (PDEs) using potential
symmetries and otherwise \cite{bl1,bl2,kb}; the use of algebraic
computing to extend Lie's methods to general third- and fourth-order
ODEs \cite{im3,ims,np}; conditional linearizability \cite{mq1,mq2};
geometric methods for symmetry analysis \cite{aa,fmq,aq} and their
use for linearization \cite{mq3,mq4}; and the development of complex
symmetry analysis \cite{sa,amq1,amq2}. Here we use the last two
developments to provide linearization for a class of systems of two
PDEs and an alternative method for linearizing a class of systems of
two ODEs that is not equivalent to the earlier method for the same
purpose \cite{amq2}.

The geometric approach for linearization of a system of two ODEs
requires that the system be at most cubic in its first derivatives
and satisfy a generalized set of Lie-Tress\'e invariant conditions
which are written in terms of the coefficients of the system of
equations \cite{mq3,mq4}. This requirement comes from regarding the
system as a projection of the geodesic equations in a flat space (in
curvilinear coordinates). The linearizing transformation is then
obtained by converting the metric in the given coordinates to one in
Cartesian coordinates. This enables us to write down the solution of
the nonlinear equations directly. In complex symmetry analysis one
studies the relationship between the algebraic properties of the
complex differential system and the corresponding real system of
ordinary or partial differential equations, including the
Cauchy-Riemann equations, which arises from the complex system via
complex splitting of the dependent and independent variables
\cite{sa,amq1,amq2}. We present linearization criteria, both
invariant and algebraic, for systems of two real ODEs and PDEs that
arise from the linearization of scalar complex second-order ODEs.
This also provides an answer to the inverse problem of those systems
of ODEs and PDEs that arise from the linearization of scalar complex
second-order ODEs. Furthermore, we provide examples of mechanical
systems that are linearizable by our procedure.

The plan of the paper is as follows. In the next Section the salient
points of complex symmetry analysis and the geometric methods
developed are mentioned. In Section 3 invariant linearizability
criteria for systems of two PDEs are given and some physical
examples provided. In the subsequent section restricted complex
transformations is introduced and used to yield invariant
linearizability criteria for systems of two ODEs. Again,
illustrative examples of mechanical systems are provided. Section 5
consists of a discussion and conclusion.

\section{Preliminaries}

A scalar ODE for a complex function of a complex variable can be
written in terms of the real and imaginary parts of the dependent
variable as functions of the real and imaginary parts of the
independent variable. As such it yields a pair of PDEs. Of course we
need to include the Cauchy-Riemann equations (C-REs), which
guarantee the differentiability of the dependent complex variable in
our system. As such we get a system of four PDEs, two of which are
first-order equations. The symmetries of the original ODE are not
identical with those of the system of PDEs. Nevertheless the
solutions of the ODE obtained by using its symmetries give the
solution of the PDE. As stated this may seem trivial. However, one
could pick a system of PDEs and check if it corresponds to an ODE.
In the present note we also specify those systems of PDEs that could
be transformed to such complex differential equations. Then one can
use the solution of the ODE to write a solution of the system of
PDEs. The correspondence between the symmetries of the system of
PDEs and the ODE is the subject of ``complex symmetry analysis"
\cite{sa,amq1}. Thus, for example, the single infinitesimal
generator for a function of one variable, ${\bf Z}=\partial/\partial
z$ becomes the pair of generators ${\bf X}=\partial/\partial x$ and
${\bf Y}=\partial/\partial y$, where $z$ is the complex independent
variable, $x$ is its real part and $y$ its imaginary part. Again the
scaling symmetry for the independent complex variable gives the
scaling symmetry for the two real variables together along with the
rotation generator in the 2-dimensional space. Similar remarks apply
for the generators involving the dependent complex variable or a
mixture of the two.

For a linear scalar ODE the conversion to the complex form and
thence to the system of two PDEs is trivial. The only complication
is that here the system becomes a set of four equations, the two
C-REs and the two PDEs for the two real variables corresponding to
the single complex variable. This can change the symmetry structure
but does not make any other substantive difference. For a nonlinear
ODE the situation changes. Now there can be non trivial mixing
between the real and imaginary parts of the complex dependent
variable so that the system of two PDEs becomes significantly
coupled. The question can arise whether the C-REs continue to hold
under some transformation of the independent and dependent
variables. This is of relevance for us as we need to use linearizing
transformations. Suppose the ODE is written for a complex analytic
function of a single variable, $u(z)$. Since the linearizing point
transformation $\mathcal{L}:(z,u) \rightarrow (Z,U)$ is analytic (by
definition), the real linearizing transformation $\mathcal{RL}:
(x,y,f,g) \rightarrow (X,Y,F,G)$ satisfies the C-REs. Note that when
we go to deal with systems of ODEs, the transformed variables are
{\it not} guaranteed to satisfy the C-REs rather they satisfy a
partial analytic structure. In that case we have to use the C-REs in
the original variables. This creates complications in the
linearization of those ODEs. Apart from this complexity the solution
of nonlinear system is still attained through a procedure analogous
to analytic continuation.

The system of geodesic equations is of second-order and
quadratically semi linear in the first derivatives and has no other
terms in it. It inherits the isometries, but can have many other
symmetries. It was noted \cite{aa} that projection of this system,
using the translational invariance symmetry of the geodesic
parameter yields a system of cubically semi linear ODEs. The
linearizability of a quadratically semi linear system of ODEs of
geodesic type is provided by regarding the coefficients of the
quadratic terms as Christoffel symbols and verifying whether the
resulting Riemann tensor is zero or not \cite{mq3}. If the system is
linearizable, one can find the linearizing transformation by taking
the coordinate transformation from the metric tensor constructed
from the Christoffel symbols \cite{fmmq} to the metric tensor in
Cartesian coordinates. The same procedure can be extended to the
projected system of cubically semi linear ODEs \cite{mq4}. The power
of this method is apparent from the fact that one can not only find
the linearizing transformations but also be able to write down the
solution of the nonlinear equation.

\section{Invariant Linearizability Criteria for Systems of PDEs}

First-order scalar ODEs are always linearizable. We only need to
deal with second-order systems for checking linearizability. The
geometric linearization of systems of two ODEs \cite{mq3,mq4} gives
a requirement that the system be (at most) cubically semi linear in
the dependent variables and satisfy a generalization of the Lie
conditions coming from the flatness of the space in which the
geodesics lie. However, this procedure does not cover all the
classes of linearizable systems as it requires a $15-$dimensional
symmetry algebra, whereas there can also be $5-, 6-, 7-$ and $8-$dimensional
algebras of linearizable systems \cite{wm}. Since the Lie symmetry
algebra for linear PDEs is infinite dimensional in general, the question
arises whether we can expect the PDEs to be (at most) cubically
semi linear as well. The geometric approach relied on the connection
between geometry and systems of second-order ODEs via the system of
geodesic equations \cite{aa,fmq}. As there is no apparent
geometrical way of extending the Lie conditions to systems of PDEs,
it is not clear what the analogues of the linearizability criteria
would be.

We extend Lie's linearizability criteria to a class of systems of
two PDEs obtainable from complex scalar ODEs. The real
transformations for linearization of this system of nonlinear PDEs
can be obtained by decomposing the complex transformations that
linearize the complex ODE. We consider a simple equation to illustrate this.
For example, we exploit the real transformation,
\begin{equation}
F=\frac{f}{f^{2}+g^{2}}-x,~G=\frac{-g}{f^{2}+g^{2}}-y,
\end{equation}%
to map a nonlinear coupled first-order system of PDEs,
\begin{eqnarray}
f_{x}+g_{y} =-2f^{2}+2g^{2},  \nonumber \\
g_{x}-f_{y} =-4fg,
\end{eqnarray}%
into a linear system
\begin{eqnarray}
F_{x}+G_{y} =0,  \nonumber \\
G_{x}-F_{y} =0.
\end{eqnarray}%
Note that the real transformation $(1)$ is equivalent to a complex
transformation%
\begin{equation}
U=(1/u)-z.
\end{equation}%
Therefore we can employ complex variables to map systems
of nonlinear PDEs into their simpler forms. It is important to mention
that the linearization above is different from the classical way of
linearizing equations because the real transformation (1) is obtained
from the complex Lie transformation (4) that may or may not be a real
Lie transformation. Therefore complex linearizability works in a different way.
A natural question following the above discussion would be ``which systems of
PDEs can be treated via complex linearization?". We provide a class of systems
of second-order PDEs which upon satisfying a set of four conditions is subject
to complex linearizability in Theorem 1. It further provides a partial answer
to another important question, i.e. ``which systems of PDEs correspond to complex
differential equations?" It is indispensable to use complex linearization in
several cases because it plays a crucial role in extracting solutions of
systems of differential equations that would have been difficult otherwise.
We highlight this feature of complex variables in the examples. In the next
Section we extend the above characteristic of complex transformations
in dealing with systems of ODEs. Namely it is shown that, if the above
transformation is restricted on a single line, then a nonlinear system of
first-order ODEs can also be linearized. In general restricted complex
transformations can be brought into play in the linearization of systems of ODEs.

For the above purpose we first prove a few results that are used
to obtain linearizability of this class of systems of PDEs.
\newline
{\bf Theorem 1}. {\it The class of systems of two second-order
partial differential equations}
\begin{eqnarray}
f_{xx}-f_{yy}+2g_{xy}=4w_1(x,y,f,g,h,l),  \nonumber \\
g_{xx}-g_{yy}-2f_{xy}=4w_2(x,y,f,g,h,l),
\end{eqnarray}
{\it where}
\begin{equation}
2h=f_x+g_y,~2l=g_x-f_y,
\end{equation}
{\it is complex-linearizable if and only if the functions}, $w_1$ {\it and}
$w_2$, {\it are at most cubic in} $h$ {\it and} $l$, {\it i.e.}
\begin{eqnarray}
f_{xx}-f_{yy}+2g_{xy}=4A^1h^3-12A^1hl^2-12A^2h^2l+4A^2l^3+4B^1h^2-
4B^1l^2-
\nonumber \\
8B^2hl+4C^1h-4C^2l+4D^1,  \nonumber \\
g_{xx}-g_{yy}-2f_{xy}=12A^1h^2l-4A^1l^3+4A^2h^3-12A^2hl^2+8B^1hl+
4B^2h^2-
\nonumber \\
4B^2l^2+4C^2h+4C^1l+4D^2
\end{eqnarray}
{\it together with the constraints on the coefficients}
\begin{eqnarray}
3A_{xx}^{1}-3A_{yy}^{1}+6A_{xy}^{2}+6C^{1}A_{x}^{1}+6C^{1}A_{y}^{2}-6A_{x}^{2}C^{2}+6C^{2}A_{y}^{1}-
\nonumber\\
6A_{f}^{1}D^{1}-6D^{1}A_{g}^{2}+6D^{2}A_{f}^{2}-6D^{2}A_{g}^{1}+6A^{1}C_{x}^{1}+6A^{1}C_{y}^{2}-6A^{2}C_{y}^{2}+
\nonumber\\
6A^{2}C_{x}^{1}+C_{ff}^{1}-C_{gg}^{1}+2C_{fg}^{2}-12A^{1}D_{f}^{1}-12A^{1}D_{g}^{2}+12A^{2}D_{f}^{2}-12A^{2}D_{g}^{1}+
\nonumber\\
2B^{1}C_{f}^{1}+2B^{1}C_{g}^{2}-2B^{2}C_{f}^{2}+2B^{2}C_{g}^{1}-4B^{1}B_{x}^{1}-4B^{1}B_{y}^{2}+4B^{2}B_{x}^{2}-4B^{2}B_{y}^{1}-
\nonumber\\
2B_{xf}^{1}-2B_{yf}^{2}-2B_{xg}^{2}+2B_{yg}^{1}=0,\nonumber \\
3A_{xx}^{2}-3A_{yy}^{2}-6A_{xy}^{1}+6C^{2}A_{x}^{1}+6C^{2}A_{y}^{2}+6A_{x}^{2}C^{1}-6C^{1}A_{y}^{1}-
\nonumber\\
6D^{2}A_{f}^{1}-6D^{2}A_{g}^{2}-6D^{1}A_{f}^{2}+6D^{1}A_{g}^{1}+6A^{2}C_{x}^{1}+6A^{2}C_{y}^{2}+6A^{1}C_{y}^{2}-
\nonumber\\
6A^{1}C_{x}^{1}+C_{ff}^{2}-C_{gg}^{2}-2C_{fg}^{1}-12A^{2}D_{f}^{1}-12A^{2}D_{g}^{2}-12A^{1}D_{f}^{2}+12A^{1}D_{g}^{1}+
\nonumber\\
2B^{2}C_{f}^{1}+2B^{2}C_{g}^{2}+2B^{1}C_{f}^{2}-2B^{1}C_{g}^{1}-4B^{2}B_{x}^{1}-4B^{2}B_{y}^{2}-4B^{1}B_{x}^{2}+
\nonumber\\
4B^{1}B_{y}^{1}-2B_{xf}^{2}+2B_{yf}^{1}+2B_{xg}^{1}-2B_{yg}^{2}=0,\nonumber \\
12D^{1}A_{x}^{1}+12D^{1}A_{y}^{2}-12D^{2}A_{x}^{2}+12D^{2}A_{y}^{1}-6D^{1}B_{f}^{1}-6D^{1}B_{g}^{2}+
\nonumber\\
6D^{2}B_{f}^{2}-6D^{2}B_{g}^{1}+6A^{1}D_{x}^{1}+6A^{1}D_{y}^{2}-6A^{2}D_{x}^{2}+6A^{2}D_{y}^{1}+
\nonumber\\
B_{xx}^{1}-B_{yy}^{1}+2B_{xy}^{2}-2C_{xf}^{1}-2C_{yf}^{2}-2C_{xg}^{2}+2C_{yg}^{1}-6B^{1}D_{f}^{1}-
\nonumber\\
6B^{1}D_{g}^{2}+6B^{2}D_{f}^{2}-6B^{2}D_{g}^{1}+3D_{ff}^{1}-3D_{gg}^{1}+6D_{fg}^{2}+4C^{1}C_{f}^{1}+
\nonumber\\
4C^{1}C_{g}^{2}-4C^{2}C_{f}^{2}+4C^{2}C_{g}^{1}-2C^{1}B_{x}^{1}-2C^{1}B_{y}^{2}+2C^{2}B_{x}^{2}-2C^{2}B_{y}^{1}=0,\nonumber \\
12D^{2}A_{x}^{1}+12D^{2}A_{y}^{2}+12D^{1}A_{x}^{2}-12D^{1}A_{y}^{1}-6D^{2}B_{f}^{1}-6D^{2}B_{g}^{2}-
\nonumber\\
6D^{1}B_{f}^{2}+6D^{1}B_{g}^{1}+6A^{2}D_{x}^{1}+6A^{2}D_{y}^{2}+6A^{1}D_{x}^{2}-6A^{1}D_{y}^{1}+
\nonumber\\
B_{xx}^{2}-B_{yy}^{2}-2B_{xy}^{1}-2C_{xf}^{2}+2C_{yf}^{1}+2C_{xg}^{1}+2C_{yg}^{2}-6B^{2}D_{f}^{1}-
\nonumber\\
6B^{2}D_{g}^{2}-6B^{1}D_{f}^{2}+6B^{1}D_{g}^{1}+3D_{ff}^{2}-3D_{gg}^{2}-6D_{fg}^{1}+4C^{2}C_{f}^{1}+
\nonumber\\
4C^{2}C_{g}^{2}+4C^{1}C_{f}^{2}-4C^{1}C_{g}^{1}-2C^{2}B_{x}^{1}-2C^{2}B_{y}^{2}-2C^{1}B_{x}^{2}+2C^{1}B_{y}^{1}=0,
\end{eqnarray}
{\it where all the coefficients} $A^{i},B^{i},C^{i}$ {\it and}
$D^{i}$,  $(i=1,2)$, {\it are functions of} $x,y,f$ and $g$.
\newline
{\bf Proof.} We firstly assume an analytic structure on the
manifold and that there exists a complex transformation that
project the functions $f(x,y)$ and $g(x,y)$ to a single complex function
$u$ of complex variable $z.$ Furthermore assume that there exists four
complex functions, $A, B, C$ and $D$, such that
\begin{eqnarray}
A(x,u)=A^1(x,y,f,g)+iA^2(x,y,f,g),  \nonumber \\
B(x,u)=B^1(x,y,f,g)+iB^2(x,y,f,g),  \nonumber \\
C(x,u)=C^1(x,y,f,g)+iC^2(x,y,f,g),  \nonumber \\
D(x,u)=D^1(x,y,f,g)+iD^2(x,y,f,g).
\end{eqnarray}
By invoking (9) we can map the system (7) to a second-order
complex differential equation
\begin{equation}
u^{\prime \prime}(x)=A(x,u)u^{\prime 3}+B(x,u)u^{\prime 2}+
C(x,u)u^{\prime}+D(x,u),
\end{equation}
which is at most cubic in its first derivative. Therefore it
satisfies Lie's linearizability criteria. Moreover the
conditions (8) can be projected down to a set of two equations
\begin{eqnarray}
3A_{zz}+3A_zC+3AC_z-3A_uD+C_{uu}-6AD_u+BC_u-2BB_z-2B_{zu}=0,
\nonumber \\
6A_zD-3B_uD+3AD_z+B_{zz}-2C_{zu}-3BD_u+3D_{uu}+2CC_u-CB_z=0;
\end{eqnarray}
which may be recognized as the Lie compatibility conditions. As the
ODE (10) is linearizable, so the system of PDEs (7) is also
linearizable.
\newline
{\bf Theorem 2.} {\it If the system of PDEs} (5) {\it admits four
real symmetries }$\mathbf{X}_{1}$, $\mathbf{Y}_{1}$, $\mathbf{X}_{2}$ {\it and} $\mathbf{Y}_{2}$, {\it
such that}
\begin{equation}
\mathbf{X}_1=\rho_1\mathbf{X}_2-\rho_2\mathbf{Y}_2,~\mathbf{Y}_1\
=\rho_1\mathbf{Y}_2+\rho_2\mathbf{X}_2,
\end{equation}
{\it for nonconstant} $\rho_1$\ {\it and} $\rho_2$, {\it and their
commutators satisfy}
\begin{equation}
[\mathbf{X}_1,\mathbf{X}_2]-[\mathbf{Y}_1,\mathbf{Y}_2]=0,
~[\mathbf{X}_1,\mathbf{Y}_2]+[\mathbf{Y}_1,\mathbf{X}_2]=0,
\end{equation}
{\it then, there exists a point transformation }$(x,y,f,g)
\longrightarrow (X,Y,F,G),$ {\it which reduces} $\mathbf{X}_{1}$, $\mathbf{Y}_{1}$,$\mathbf{X}_{2}$ {\it and} $\mathbf{Y}_{2}$ {\it to their canonical form}
\begin{equation}
\mathbf{X}_1=\frac{\partial}{\partial F},\quad \mathbf{Y}_1=\frac{\partial}
{\partial G},\quad \mathbf{X}_2=X\frac{\partial}{\partial F}+Y\frac{\partial}
{\partial G},\quad \mathbf{Y}_2=Y\frac{\partial}{\partial F}-X\frac{\partial}
{\partial G},
\end{equation}
{\it and the system} (5) {\it can be reduced to the linear form}
\begin{eqnarray}
F_{XX}-F_{YY}+2G_{XY}=4W_1(X,Y),  \nonumber \\
G_{XX}-G_{YY}-2F_{XY}=4W_2(X,Y).
\end{eqnarray}
{\bf Proof.} Suppose that $\mathbf{X}_a+i\mathbf{Y}_a=$ $\mathbf{Z}
_a,$ for $a=1,2$. Then equation (13) can be replaced by
$[\mathbf{Z}_1,\mathbf{Z}_2]=0,$ which implies that the two complex
symmetries $\mathbf{Z}_1,~\mathbf{Z}_2$, commute with each other.
Further, setting $\mathbf{Z}_1=\rho(z,u)\mathbf{Z}_2$ for a
nonconstant complex function $\rho$, justifies equation (12). It is
proved by Lie (see \cite{fmm}) that every scalar second-order ODE (10) that
admits two commuting symmetries such that
$\mathbf{Z}_{1}=\rho (z,u) \mathbf{Z}_{2}$ can be transformed into a
linear ODE $U''=W(\zeta)$, by applying the point
transformation, $\zeta=\zeta(z,u),U=U(z,u)$ which reduces
$\mathbf{Z}_1$ and $\mathbf{Z}_2$ to their canonical forms
\begin{equation}
\mathbf{Z}_1=\frac{\partial}{\partial U},~\mathbf{Z}_2=\zeta
\frac{\partial}{\partial U}.
\end{equation}
The point transformation, $(x,y,f,g)\longrightarrow (X,Y,F,G),$ can be
obtained for complex transformation which then can be used to convert
system (5) into the linear form (15).
\newline
{\bf Examples}
\newline
{\bf 1.} \textit{Higher-dimensional Coupled System of Modified
Lane-Emden Type:}\newline The Lane-Emden equation arises in the
study of stellar structures \cite{chand}. We first investigate the
linearizability of a system of two cubically semi linear PDEs
\begin{eqnarray}
f_{xx}-f_{yy}+2g_{xy}=-12fh+12gl-4f^3+12fg^2,\quad f_{x}=g_{y},  \nonumber \\
g_{xx}-g_{yy}-2f_{xy}=-12gh-12fl-12f^2g+4g^3,\quad f_{y}=-g_{x}.
\end{eqnarray}
The above system can be regarded as a special case of a
higher-dimensional Lane-Emden system equipped with an analytic structure.
The coefficients satisfy linearizability conditions (8). Therefore
the above system is linearizable. Thus one can check the
linearizability of a class of systems of PDEs analogously to Lie's
technique. In order to construct the transformation we employ the
symmetries of (17). It admits Lie symmetries
\begin{eqnarray}
\mathbf{X}_1=\frac{\partial}{\partial x},\quad \mathbf{Y}_1=
\frac{\partial}{\partial y},  \nonumber \\
\mathbf{X}_2=x\frac{\partial}{\partial x}+
y\frac{\partial}{\partial y}-f\frac{\partial}{\partial f}-
g\frac{\partial}{\partial g},\quad
\mathbf{Y}_2=y\frac{\partial}{\partial x}-
x\frac{\partial}{\partial y}-g\frac{\partial}{\partial f}+
f\frac{\partial }{\partial g},
\end{eqnarray}
which satisfy all the conditions of Theorem 2. Therefore (18)
leads to the linearizing transformation,
\begin{eqnarray}
X=x-\frac{f}{f^2+g^2},\quad Y=y+\frac{g}{f^2+g^2}, \nonumber \\
F=\frac{1}{2}(x^2-y^2)-\frac{f^2+g^2}{xf+yg},\quad
G=xy-\frac{f^2+g^2}{yf-xg},
\end{eqnarray}
and we deduce the following system of linear PDEs
\begin{eqnarray}
F_{XX}-F_{YY}+2G_{XY}=0,  \nonumber\\
G_{XX}-G_{YY}-2F_{XY}=0.
\end{eqnarray}
{\bf 2.} Consider the nonlinear anharmonic oscillator system given
by the PDEs
\begin{eqnarray}
f(f_{xx}-f_{yy}+2g_{xy})-g(g_{xx}-g_{yy}-2f_{xy})
=4(h^2-l^2)-4(f^2-g^2)w_1+8fgw_2,  \nonumber \\
f(g_{xx}-g_{yy}-2f_{xy})+g(f_{xx}-f_{yy}+2g_{xy})
=8hl-8fgw_1-4(f^2-g^2)w_2, \nonumber\\
f_{x}=g_{y},\quad f_{y}=-g_{x},\nonumber \\
\end{eqnarray}
where both $w_1$ and $w_2$ are arbitrary functions of $x$ and $y$.
We can transform the above system into a system of the form (7).
Consequently it can be verified that the coefficients of $h$ and $l$
satisfy the conditions (8). Therefore the above system is
linearizable. Notice that we can convert the system (21) into the
second-order complex nonlinear anharmonic oscillator ODE
\begin{equation}
uu^{\prime \prime }+w(z)u^2=u^{\prime 2}.
\end{equation}
It admits the symmetries
\begin{equation}
\mathbf{Z}_1=zu\frac{\partial}{\partial u},\quad \mathbf{Z}_2=
u\frac{\partial}{\partial u},
\end{equation}
which yields the complex transformation
\begin{equation}
Z=\frac{1}{z}, \quad U=\frac{1}{z}\log u.
\end{equation}
The above transformation gives the real transformation
\begin{eqnarray}
X=\frac{x}{x^2+y^2},\quad Y=\frac{-y}{x^2+y^2},  \nonumber \\
F=\frac{1/2}{x^2+y^2}\left (x\ln \left (f^2+g^2 \right )+2
y\arctan\left(g/f\right )\right ),  \nonumber \\
G=\frac{1/2}{x^2+y^2}\left (2x\arctan\left (g/f \right )-
y\ln\left (f^2+g^2\right )\right ),
\end{eqnarray}
that transforms system (21) into the linear form
\begin{eqnarray}
F_{XX}-F_{YY}+2G_{XY}=-\frac{4}{X^2+Y^2}\left ((X^3-3XY^2)w_1-
(Y^3-3X^2Y)w_2\right ),\nonumber \\
G_{XX}-G_{YY}-2F_{XY}=-\frac{4}{X^2+Y^2}\left ((X^3-3XY^2)w_2+
(Y^3-3X^2Y)w_1\right),
\end{eqnarray}
where
\begin{eqnarray}
w_{1}=w_{1}\left(\frac{X}{X^{2}+Y^{2}},\frac{-Y}{X^{2}+Y^{2}}\right), \nonumber\\
w_{2}=w_{2}\left(\frac{X}{X^{2}+Y^{2}},\frac{-Y}{X^{2}+Y^{2}}\right ).
\end{eqnarray}
\textbf{3.} We use the transformation
\begin{equation}
X =2f-x^{2}+y^{2},~Y=2g-2xy,~F=x,~G=y,
\end{equation}
to linearize the system of PDEs
\begin{eqnarray}
f_{xx}-f_{yy}+2g_{xy}=4+4\left((h-x)^2-(l-y)^2\right )w_1-8(h-x)(l-y)w_2,\nonumber \\
g_{xx}-g_{yy}-2f_{xy}=8(h-x)(l-y)w_1+4\left((h-x)^2-(l-y)^2\right)w_2,\nonumber \\
f_{x}=g_{y},\quad f_{y}=-g_{x},
\end{eqnarray}
where $w_1$ and $w_2$ are arbitrary functions of the mixed variables
$2f-x^{2}+y^{2}$ and $2g-2xy$. The above system possesses the
symmetries
\begin{eqnarray}
\mathbf{X}_1=\frac{\partial}{\partial x}+x\frac{\partial}
{\partial f}+y\frac{\partial}{\partial g},\quad \mathbf{Y}_1=
x\frac{\partial}{\partial g}-y\frac{\partial}{\partial f}
-\frac{\partial}{\partial y}, \quad \nonumber \\
\mathbf{X}_2=x\frac{\partial}{\partial x}+y\frac{\partial}
{\partial y}+(x^2-y^2)\frac{\partial}{\partial f}+2xy\frac{\partial}
{\partial g}, \nonumber \\
\mathbf{Y}_2=y\frac{\partial}{\partial x}-x\frac{\partial}
{\partial y}+2xy\frac{\partial}{\partial f}-(x^{2}-y^{2})
\frac{\partial}{\partial g}.
\end{eqnarray}
The linearized system of PDEs is
\begin{eqnarray}
F_{XX}-F_{YY}+2G_{XY}=-2(Hw_1-Lw_2),  \nonumber \\
G_{XX}-G_{YY}-2F_{XY}=-2(Hw_2+Lw_1).
\end{eqnarray}
These conclude our examples. In the next Section we present another important use of the complex
method.
\section{Restricted Complex Transformations and Invariant \\
Linearizability Criteria for Systems of Two ODEs}

In the previous section we utilized complex functions of complex
variables to obtain linearizing transformations for systems of PDEs.
If we restrict our complex functions to depend upon a single real
variable, then it would generate different transformations. They
transform systems of ODEs into other systems of ODEs. We extend our
three-dimensional space of two dependent and one independent
variables to a two-complex-dimensional space. In the intermediate
steps we move off the real line to obtain our linearizing
transformations. Then the solutions of systems of ODEs are recovered
by restricting the independent variable to the real line. The
procedure is reminiscent of analytic continuation. In the process we
loose the C-REs so our system is not the same as that with which we
started. It is very interesting to see the invariance of analytic
structure under the restricted complex transformation. In this case
the ODE is written for an analytic function of a single real
variable, $u(x)$. It is important to understand how the complex
transformation works in the restricted domain. The linearizing
complex point transformation $\mathcal{L}:(x,u) \rightarrow
(\chi,U)$ can exhibit dual nature in terms of its analyticity which
yields the following real linearizing transformation $\mathcal{RL}:
(x,f,g) \rightarrow (\chi,\Upsilon,\zeta)$. Because in the
transformed variables $\chi$ can be either complex or real. If it is
real, then the transformed variables satisfy C-REs only in $f$ and
$g$ in which case there is a partial analytic structure on the
transformed variables, but, if $\chi$ is complex, then the complete
analytic structure is restored on the transformed variables. We
illustrate this important feature in the two examples below. Thus we
cannot guarantee that at the end a similar linearized system can be
obtained via the linearizing transformations from other approaches
even though we do get the solution by linearization (in the
complex).

The invariance properties of a system of ODEs,
\begin{eqnarray}
f''=w_1(x,f,g,f',g'), \nonumber \\
g''=w_2(x,f,g,f',g'),
\end{eqnarray}
have been investigated by using complex functions in \cite{sa,amq1}.
The idea is to make use of the transformation,
\begin{equation}
u(x)=f(x)+ig(x),~~w(x,u)=w_{1}(x,f,g)+iw_{2}(x,f,g),
\end{equation}
to convert the system (32) into the single ODE,
\begin{equation}
u''(x)=w(x,u,u'),
\end{equation}
and then use the standard Lie procedure of linearization. We call
this procedure {\it complex linearization} even though we have no
guarantee that the transformed system can be linearized via other
approaches. To comprehend the connection between complex
transformations and linearization we firstly revisit an earlier
example from the previous section.

Consider a first-order two-dimensional Riccati system%
\begin{eqnarray}
f^{\prime }=-f^{2}+g^{2}, \nonumber\\
g^{\prime }=-2fg,
\end{eqnarray}%
which is a special case of a general Riccati system in two
dimensions. A natural question arises: which transformation can
linearize the above system? If such a transformation exists, then
how can we find it. The use of the transformation,
\begin{equation}
\Upsilon =\frac{f}{f^{2}+g^{2}}-x,\quad \zeta =\frac{-g}{f^{2}+g^{2}},
\end{equation}%
maps system $(35)$ into the simplest system%
\begin{equation}
\Upsilon ^{\prime }=0,\quad\zeta ^{\prime }=0.
\end{equation}%
The transformation $(36)$ \textit{is} a mere consequence of the same
complex transformation $(4)$ that we used to linearize system of PDEs $%
(2)$ in the remaining Section. The only difference is that we restricted the
dependent variable to a single real line. This indicates to us a significant use of restricted
complex transformations in converting systems of nonlinear
equations into their linear analogues. Theorem 3 provides a class of systems of
two second-order ODEs that can be dealt via complex variables. We now present the basic theorem for complex linearization.
\newline
{\bf Theorem 3.} {\it The necessary and sufficient condition for a
system of two second-order ODEs of the form}
\begin{eqnarray}
f''=A^1f'^3-3A^1f'g'^2-3A^2f'^2g'+A^2g'^3+B^1f'^2-B^1g'^2-2B^2f'g'+\nonumber \\
C^1f'-C^2g'+D^1,  \nonumber \\
g''=3A^1f'^2g'-A^1g'^3+A^2f'^3-3A^2f'g'^2+2B^1f'g'+B^2f'^2-B^2g'^2+\nonumber \\
C^2f'+C^1g'+D^2,
\end{eqnarray}
{\it where} $A^{i}$, $B^{i}$, $C^{i}$, $D^{i}$, $(i=1,2)$
{\it are functions of the variables} $x,y,f$ {\it and} $g$, {\it to be
solvable by complex linearization is that the coefficients satisfy
the conditions}
\begin{eqnarray}
12A_{xx}^{1}+12C^{1}A_{x}^{1}-12A_{x}^{2}C^{2}-6A_{f}^{1}D^{1}-6D^{1}A_{g}^{2}+6D^{2}A_{f}^{2}-6D^{2}A_{g}^{1}+
\nonumber \\
12A^{1}C_{x}^{1}-12A^{2}C_{x}^{2}+C_{ff}^{1}-C_{gg}^{1}+2C_{fg}^{2}-12A^{1}D_{f}^{1}-12A^{1}D_{g}^{2}+
\nonumber \\
12A^{2}D_{f}^{2}-12A^{2}D_{g}^{1}+2B^{1}C_{f}^{1}+2B^{1}C_{g}^{2}-2B^{2}C_{f}^{2}+2B^{2}C_{g}^{1}-8B^{1}B_{x}^{1}+
\nonumber\\
8B^{2}B_{x}^{2}-4B_{xf}^{1}-4B_{xg}^{2}=0,
\nonumber \\
12A_{xx}^{2}+12C^{2}A_{x}^{1}+12A_{x}^{2}C^{1}-6D^{2}A_{f}^{1}-6D^{2}A_{g}^{2}-6D^{1}A_{f}^{2}+6D^{1}A_{g}^{1}+
\nonumber \\
12A^{2}C_{x}^{1}+12A^{1}C_{x}^{2}+C_{ff}^{2}-C_{gg}^{2}-2C_{fg}^{1}-12A^{2}D_{f}^{1}-12A^{2}D_{g}^{2}-
\nonumber \\
12A^{1}D_{f}^{2}+12A^{1}D_{g}^{1}+2B^{2}C_{f}^{1}+2B^{2}C_{g}^{2}+2B^{1}C_{f}^{2}-2B^{1}C_{g}^{1}-8B^{2}B_{x}^{1}-
\nonumber\\
8B^{1}B_{x}^{2}-4B_{xf}^{2}+4B_{xg}^{1}=0,
\nonumber \\
24D^{1}A_{x}^{1}-24D^{2}A_{x}^{2}-6D^{1}B_{f}^{1}-6D^{1}B_{g}^{2}+6D^{2}B_{f}^{2}-6D^{2}B_{g}^{1}+
\nonumber \\
12A^{1}D_{x}^{1}-12A^{2}D_{x}^{2}+4B_{xx}^{1}-4C_{xf}^{1}-4C_{xg}^{2}-6B^{1}D_{f}^{1}-6B^{1}D_{g}^{2}+
\nonumber \\
6B^{2}D_{g}^{2}-6B^{2}D_{g}^{1}+3D_{ff}^{1}-3D_{gg}^{1}+6D_{fg}^{2}+4C^{1}C_{f}^{1}+4C^{1}C_{g}^{2}-
\nonumber\\
4C^{2}C_{f}^{2}+4C^{2}C_{g}^{1}-4C^{1}B_{x}^{1}+4C^{2}B_{x}^{2}=0,
\nonumber \\
24D^{2}A_{x}^{1}+24D^{1}A_{x}^{2}-6D^{2}B_{f}^{1}-6D^{2}B_{g}^{2}-6D^{1}B_{f}^{2}+6D^{1}B_{g}^{1}+
\nonumber \\
12A^{2}D_{x}^{1}+12A^{1}D_{x}^{2}+4B_{xx}^{2}-4C_{xf}^{2}+4C_{xg}^{1}-6B^{2}D_{f}^{1}-6B^{2}D_{g}^{2}-
\nonumber \\
6B^{1}D_{f}^{2}+6B^{1}D_{g}^{1}+3D_{ff}^{2}-3D_{gg}^{2}-6D_{fg}^{1}+4C^{2}C_{f}^{1}-4C^{2}C_{g}^{2}+
\nonumber\\
4C^{1}C_{f}^{2}-4C^{1}C_{g}^{1}-4C^{2}B_{x}^{1}-4C^{1}B_{x}^{2}=0.
\end{eqnarray}
\textbf{Proof.} Suppose that there exists complex functions,
\begin{eqnarray}
A(x,u)=A^{1}(x,f,g)+iA^{2}(x,f,g), \nonumber \\
B(x,u)=B^{1}(x,f,g)+iB^{2}(x,f,g), \nonumber \\
C(x,u)=C^{1}(x,f,g)+iC^{2}(x,f,g), \nonumber \\
D(x,u)=D^{1}(x,f,g)+iD^{2}(x,f,g),
\end{eqnarray}
such that the above system can be mapped into the second-order ODE
\begin{equation}
u''(x)=A(x,u)u'^3+B(x,u)u'^2+C(x,u)u'+D(x,u),
\end{equation}
which is at most cubic in $u'$ and therefore satisfies the necessary
condition of linearizability. To check the sufficient conditions
the set of equations (39) is transformed into the equations
\begin{eqnarray}
3A_{xx}+3A_{x}C+3AC_{x}-3A_{u}D+C_{uu}-6AD_{u}+BC_{u}-2BB_{x}-
2B_{xu}=0, \nonumber \\
6A_{x}D-3B_{u}D+3AD_{x}+B_{xx}-2C_{xu}-3BD_{u}+3D_{{u}u}+2CC_{u}-
CB_{x}=0,
\end{eqnarray}
which are the Lie conditions. Since the ODE (41) is linearizable,
we can obtain its solution by the geometric method. This solution
can now be written as the pair of real functions, $f$ and $g$.
Hence the system of ODEs (38) can be solved by complex
linearization.

To see how the complex variable approach works we present some
illustrative examples of two dimensional systems of ODEs.
\newline
\textbf{Examples}
\newline
\textbf{1.} Consider the system of ODEs
\begin{eqnarray}
ff''-gg''=f'^2-g'^2 + (f^2-g^2) w_{1}(x) - 2fgw_2(x),  \nonumber \\
fg''+gf''=2f'g' + 2fg w_1(x) + (f^2-g^2) w_2(x),
\end{eqnarray}
where $w_1$ and $w_2$ are arbitrary functions of the time variable
$x$. This was discussed in \cite{faq} for an anharmonic oscillator
with $w_1=x^2$ and $w_2=0$. The complex symmetry analysis gave
interesting insights into the oscillator dynamics. Since the
coefficients satisfy the linearizability conditions therefore the
above system is solvable by complex linearization. In fact the
system is linearizable. The transformation for our purpose is
\begin{equation}
\chi =\frac{x}{x^{2}+f^{2}},\quad \Upsilon =\frac{1}{2x}\ln (f^{2}+g^{2}),\quad
\zeta =\frac{1}{x}\arctan \left (\frac{g}{f} \right ),
\end{equation}
and it reduces system (43) into the linear system of ODEs
\begin{equation}
\Upsilon''=\frac{1}{\chi}w_1,\quad\zeta''=\frac{1}{\chi}w_2,
\end{equation}
where
\begin{equation}
w_1\equiv w_1(1/\chi),\quad w_2\equiv w_2(1/\chi).
\end{equation}
\newline
\textbf{2.} \textit{Two-dimensional Coupled Modified Emden System:}
\newline The modified Emden equation possess numerous dynamical
properties in nonlinear oscillations. In \cite{ch1} {\it
Chandrasekar et al}. explore an important characteristic of such an
equation in which the frequency of oscillation is independent of
amplitude and remains the same as that of the linear oscillator.
They showed that the amplitude dependence of the frequency is not a
fundamental property of nonlinear dynamical phenomena. Later, in
\cite{ch2}, they extended the above results and investigated the
dynamical properties of N-coupled nonlinear oscillator of Lienard
type. Those systems which possess a Hamiltonian structure can be
transformed into systems of uncoupled harmonic oscillators via
contact transformations. We investigate the integrability of a
system of two cubically semilinear coupled ODEs of Emden type
\begin{eqnarray}
f''=-3ff'+3gg'-f^3+3fg^2, \nonumber \\
g''=-3gf'-3fg'-3f^2g+g^3.
\end{eqnarray}
The coefficients satisfy the linearizability conditions. Hence the
above system is solvable by complex linearization. This system
corresponds to the cubically semilinear modified Emden ODE
\begin{equation}
u''+3uu'+u^3=0,
\end{equation}
which is linearizable as it satisfies the Lie conditions. This
equation arises in many applications (see e.g. \cite{fmm,ml1}).
Equation (48) admits the two noncommuting complex symmetries
\begin{equation}
\mathbf{Z}_1=\frac{\partial}{\partial x},\quad \mathbf{Z}_2=
x\frac{\partial}{\partial x}-u\frac{\partial}{\partial u},
\end{equation}
which can be used to write down the transformation,
\begin{equation}
\chi =x-\frac{1}{u},\quad U=\frac{x^{2}}{2}-\frac{x}{u},
\end{equation}
that transforms (48) into the free-particle equation $U''=0$. This
transformation seems odd as $x$ is real while $\chi $ is complex. The
point is that we start and end with a complex independent variable
{\it restricted to the real line}, but in the intervening steps the
variable moves off it. The procedure is reminiscent of analytic
continuation. To check its consistency we express the solution of
(47) in the new coordinates $(\chi ,U)$
\begin{equation}
U=\alpha \chi +\beta ,
\end{equation}
where $\alpha $ and $\beta $ are complex constants. In coordinates
$(x,u)$ the above equation yields
\begin{equation}
u=\frac{2(x-\alpha)}{x^2-2\alpha x-2\beta},
\end{equation}
which satisfies $(48)$. It generates the solution of the system
$(47)$
\begin{eqnarray}
f=\frac{2(x-\alpha_1)(x^2-2\alpha_1x-2\beta_1)+4\alpha_2
(\alpha_2x+\beta_2)}{(x^2-2\alpha_1x-2\beta_1)^2+
(2\alpha_2x+2\beta_2)^2},  \nonumber \\
g=\frac{4(x-\alpha_1)(\alpha_2x+\beta_2)-2\alpha_2
(x^2-2\alpha_1x-2\beta_1)}{(x^2-2\alpha_1x-2\beta_1)^2+
(2\alpha_2x+2\beta_2)^2},
\end{eqnarray}
which would have been difficult to obtain by other means. Note that
here the C-REs are no longer preserved by the transformation.
\newline
\textbf{3.} Now consider the Newtonian system of ODEs with velocity
dependent forces
\begin{eqnarray}
f''=1+\left((f'-x)^2-g'^2\right)w_1-2(f'-x)g'w_2,  \nonumber \\
g''=2(f'-x)g'w_1+\left((f'-x)^2-g'^2\right)w_2,
\end{eqnarray}
where
\begin{equation}
w_1\equiv w_1(2f-x^2,2g), \quad w_2\equiv w_2(2f-x^2,2g).
\end{equation}
It can be verified that system $(54)$ can be solved by complex
linearization as it satisfies the criteria of Theorem 3. It
corresponds to the complex Newtonian equation
\begin{equation}
u''=1+(u'-x)^2w(2u-x^2)
\end{equation}
with quadratic velocity dependent forces. This ODE admits the
complex Lie symmetries
\begin{equation}
\mathbf{Z}_1=\frac{\partial}{\partial x}+x\frac{\partial}{\partial u},\quad
\mathbf{Z}_2=x\frac{\partial}{\partial x}+
x^2\frac{\partial}{\partial u}.
\end{equation}
Equation $(56)$ is linearized by the complex transformation
\begin{equation}
\chi =2u-x^2, \quad U=x,
\end{equation}
to become
\begin{equation}
2U''=-U'w(\chi ).
\end{equation}
Again (58) is reminiscent of analytic continuation as it is
from (real, complex) to (complex, real). To check consistency we
may take $w=1$, i.e., $w_1=1,w_2=0,$ to obtain
\begin{equation}
u=\alpha +\ln 2+\frac{x^2}{2}-\ln(\beta -x),
\end{equation}
where $\alpha $ and $\beta $ are complex constants. Putting
$w_1=1$ and $w_2=0$ in (54) we get the system
\begin{eqnarray}
f''=1+x^2+f'^2-g'^2-2xf', \nonumber \\
g^{\prime \prime}=2f'g'-xg',
\end{eqnarray}
with the general solution
\begin{eqnarray}
f=\alpha_1-\ln 2+\frac{x^2}{2}-\frac{1}{2}\ln\left((\beta_1-x)^2+\beta_2^2\right), \nonumber \\
g=\alpha_2-\arctan\left(\frac{\beta_2}{\beta_2-x}\right).
\end{eqnarray}
Also, if we take
\begin{equation}
w(2u-x^2)=\frac{1}{2u-x^2},
\end{equation}
the solution of (56) is
\begin{equation}
u=\frac{x^2}{2}+\sqrt{\frac{\beta -x}{2\alpha}},
\end{equation}
which gives the solution
\begin{eqnarray}
f(x)=\frac{x^{2}}{2}+R(x)\cos (\theta(x)),
\nonumber \\
g(x)=R(x) \sin\left (\theta(x)\right),
\end{eqnarray}
where
\begin{eqnarray}
R(x)=\sqrt{\frac{(\alpha _{1}(\beta _{1}-x)+
\beta _{2}\alpha_{2})^{2}+(\beta _{2}\alpha _{2}-\alpha _{2}
(\beta _{1}-x))^{2}}{2(\alpha _{1}^{2}+\alpha _{2}^{2})}}, \nonumber \\
\theta(x)= \frac{\beta _{2}\alpha _{2}-\alpha_{2}(\beta _{1}-x)}
{2(\alpha _{1}(\beta _{1}-x)+\beta _{2}\alpha _{2})}, \nonumber
\end{eqnarray}
of the system
\begin{eqnarray}
f''=1+\frac{(2f-x^2)((f'-x)^2-g'^2)}{(2f-x^2)^2+4g^2}+
\frac{4gg'(f'-x)}{(2f-x^2)^2+4g^2},  \nonumber \\
g''=\frac{2g'(f'-x)(2f-x^2)}{(2f-x^2)^2+4g^2}-
\frac{2g((f'-x)^2-g'^2)}{(2f-x^2)^2+4g^2}.
\end{eqnarray}
Note that $w$ is an arbitrary complex function that gives a class of
systems of ODEs that correspond to (56). Thus the linearization
of a general equation encodes the linearization of a large class of
systems of ODEs.

\section{Conclusion and Discussion}

Though the linearizaton procedure for a scalar ODE was fully provided
by Lie, there is no such complete characterization and procedure
available more generally \cite{bl1,bl2,kb}, \cite{fmm}$-$\cite{ml2}).
Geometry gives a procedure that not only provides the invariant
characterization but also the solution of systems of ODEs
\cite{mq3,mq4}. However, it only applies to the class of maximum
symmetry and it is known that some less symmetric systems of ODEs
are also linearizable \cite{wm}. A method that retains the power of
geometry, but applies to the less symmetric cases is needed.

In this paper we used complex scalar ODEs to write equivalent
systems of PDEs and ODEs and then required that the original ODEs be
linearizable. We provided various examples, mainly of mechanical
systems, to illustrate the power and use of the method presented.
For the PDEs we have a guarantee that the C-REs are preserved under
the linearizing transformation. However, for the ODEs we have no
such guarantee. It seems that, when the C-REs are preserved, we
obtain linearizing transformations and when they are not preserved
we are not able to linearize the system. However, the solution of
the system obtained by linearization of the complex scalar ODE is
still valid for the equivalent system of PDEs or ODEs.

We called the procedure for ODEs solution by complex linearization.
This is an example of Penrose's ``complex magic" \cite{pen}. The
solutions are again obtained using the geometric linearization
procedure. It is hoped that the complex procedures, augmenting the
real linearization for maximally symmetric systems, will give all
the classes of linearizable systems of ODEs. This line of
investigation is being pursued \cite{amq2,sq}.

For PDEs we get a systematic procedure to linearize the system and
obtain the solution. However, it is clear that we do not get all
possible solutions. This can be proved by considering the linear
complex ODE, $u''=0$, and writing the corresponding system of PDEs,
\begin{eqnarray}
f_{xx}-f_{yy}+2g_{xy}=0, \nonumber \\
g_{xx}-g_{yy}-2f_{xy}=0.
\end{eqnarray}
Since there are only two arbitrary (complex) constants that can
appear for the scalar ODE but infinitely many linearly independent
solutions of the system of PDEs, we see that the complex
linearization procedure cannot exhaust the solutions of the system of
PDEs.

One should be able to extend the complex linearization procedure to
third- and fourth-order systems by using the results for the
corresponding scalar ODEs \cite{im3,ims} straightforwardly. Also
the extension to conditional linearizability of systems could be
obtained \cite{mnq,mq1,mq2}. However, in the former case the power
of the geometric approach would be lost. In the latter case, if the
root equation is second-order, the geometric method would provide the
solution of the system.

A bigger problem is the extension to higher-dimensional systems. The
extension to dimensions of $2n$ can be obtained by iterative use of
the complex method. Starting with a real system of $n$-dimensions
and converting to complex variables, we could get $2n$-dimensional
systems as well. However, it would appear that the iterative
procedure and the complexification of the system would yield
different linearizable classes in general. Furthermore, this would not
provide a means of dealing with systems of odd dimensions. It would be
interesting to explore the various ramifications of extension of
complex linearization to higher dimensions.

\section*{Acknowledgments}

\qquad SA is most grateful to NUST and DECMA in providing financial
assistance for his stay at the University of the Witwatersrand, Johannesburg, South Africa, where this
work was initiated. AQ acknowledges the School of Computational and
Applied Mathematics (DECMA) for funding his stay at the university.

\end{document}